\documentclass{article}
\usepackage{amsmath}

\begin{document}

\large{
\bigskip
\begin{center}
\textbf{Investigation of one  boundary-value problem for elliptic
type equation}
\end{center}

\begin{center}
\textbf{Asadova O.H.}
\end{center}
\begin{center}
Faculty of Applied Mathematics and Cybernetics  Baku State
University,\\ Z.Khalilov str.23, Baku AZ1148, Azerbaijan,
e-mail:aahmado7@rambler.ru
\end{center}

\bigskip

As is known one of the stages of the solution of a mixed problem by the
contour integral method [1]-[2] is the solution of a spectral problem
corresponding to the given mixed problem.

If we consider the problem for a parabolic  Petrovski equation
then by means of special potentials the solution of a spectral
problem is reduced to the solution of a system of regular integral
equations with respect to unknown densities.

By virtue of the parabolic character of the equation, the kernels of these
potentials decrease well at great values of a parameter from some infinite
part of a complex plane and have weak pointwise singularity in space
variable.

The goal of the paper is to study a spectral problem corresponding to the
mixed problem for a sixth order weak parabolic equation. For simplicity we
consider an equation with constant coefficients, namely, we consider a
boundary value problem on finding of the solution of the equation

\begin{equation*}
A_{0}\Delta ^{3}u\left( x,\lambda \right) +A_{1}\lambda ^{2}\Delta
^{2}u\left( x,\lambda \right) +A_{2}\lambda ^{4}\Delta u\left( x,\lambda
\right) +\lambda ^{6}u\left( x,\lambda \right) =0,\ \eqno(1)
\end{equation*}%
satisfying the boundary conditions

\begin{equation*}
\underset{x\rightarrow z\in \tau }{\lim }u\left( x,\lambda \right) =\varphi
_{0}\left( z,\lambda \right) ;\;\underset{x\rightarrow z\in \tau }{\lim }%
\frac{du\left( x,\lambda \right) }{dn_{z}}=\varphi _{1}\left( z,\lambda
\right)
\end{equation*}

\begin{equation*}
\underset{x\rightarrow z\in \tau }{\lim }\frac{d}{dn_{z}}\Delta ^{2}u\left(
x,\lambda \right) =\varphi _{2}\left( z,\lambda \right) ,\eqno(2)
\end{equation*}%
where $x=\left( x_{1},x_{2}\right) $ is a point of domain $D$ with boundary $%
\tau ,\ A_{k}\;\left( k=\overline{0,2}\right) $ are constant numbers, $%
\lambda $ is a complex parameter, $n_{z}$ is the direction of an
inside normal to the boundary $\tau $ of the domain $D$ at the
point $z\in \tau $.

The fulfillment of the following conditions is assumed:

1) The roots of the characteristic equation

\begin{equation*}
\nu ^{3}-A_{2}\nu ^{2}+A_{1}\nu -A_{0}=0,\ \eqno(3)
\end{equation*}%
corresponding to equation (1) are such that a real part of even if one of
them equals zero, and the others are negative, i.e., $ {Re}\nu _{1}=0,\;%
 {Re}\nu _{k}<0\;\left( k=2,3\right) $

2) The boundary functions $\varphi _{s}\left( z,\lambda \right) \;\left( s=%
\overline{0,2}\right) $ have continuous derivatives with respect
to $z$ up to three
order for $z\in \tau $, with analytic by $\lambda $ in $%
R_{\delta }$ functions and converging to zero as $\left\vert \lambda
\right\vert \rightarrow \infty $.

3) The boundary $\tau $ of the domain $D$ is a Lyapunov line. We look for
the solution of problem (1), (2) in the form of sum of potentials

\begin{equation*}
u\left( x,\lambda \right) =w_{1}\left( x,\lambda \right) +w_{2}\left(
x,\lambda \right) +w_{3}\left( x,\lambda \right) ,\ \eqno(4)
\end{equation*}%
where $w_{k}\left( x,\lambda \right) \;\left( k=\overline{1,3}\right) $ are
special potentials defined by the formulae

\begin{equation*}
w_{1}\left( x,\lambda \right) =\underset{\tau }{\int }\mathcal{P}_{0}\left(
x-y;\lambda \right) \mu _{1}\left( y;\lambda \right) d\tau _{y},\;\eqno(5)
\end{equation*}

\begin{equation*}
w_{2}\left( x,\lambda \right) =\underset{\tau }{\int }\mathcal{P}_{1}\left(
x-y;\lambda \right) \mu _{2}\left( y;\lambda \right) d\tau _{y}\eqno(6)
\end{equation*}

\begin{equation*}
w_{3}\left( x,\lambda \right) =\underset{\tau }{\int }\mathcal{P}_{2}\left(
x-y;\lambda \right) \mu _{3}\left( y;\lambda \right) d\tau _{y}\eqno(7)
\end{equation*}%
where $\mathcal{P}\left( x-y;\lambda \right) $ is a fundamental solution,
and $\mathcal{P}_{s}\left( x-y;\lambda \right) \ \left( s=1,2\right) $ are
particular solutions of equation (1) (see [3], [4])

\begin{equation*}
\mathcal{P}_{0}\left( x-y;\lambda \right) =-\frac{1}{4\pi \lambda ^{4}}%
\underset{k=1}{\overset{3}{\sum }}\frac{\nu _{k}\mathcal{K}_{0}\left( \frac{%
\lambda \left\vert x-y\right\vert }{\sqrt{-\nu _{k}}}\right) }{\overset{3}{%
\underset{\underset{s\neq k}{s=1}}{\prod }}\left( \nu _{k}-\nu _{s}\right) }%
,
\end{equation*}

\bigskip
\begin{equation*}
\mathcal{P}_{1}\left( x-y;\lambda \right) =-\frac{1}{4\pi \lambda ^{4}}%
\underset{k=1}{\overset{3}{\sum }}\frac{\nu _{k}^{2}\mathcal{K}_{0}\left(
\frac{\lambda \left\vert x-y\right\vert }{\sqrt{-\nu _{k}}}\right) }{\overset%
{3}{\underset{\underset{s\neq k}{s=1}}{\prod }}\left( \nu _{k}-\nu
_{s}\right) },
\end{equation*}

\begin{equation*}
\mathcal{P}_{2}\left( x-y;\lambda \right) =\left[ \mathcal{P}_{2}\left(
x-y;\lambda \right) -\frac{2A_{0}}{3\lambda ^{2}}\frac{d^{2}}{dn_{y}^{2}}%
\mathcal{P}_{3}^{\ast }\left( x-y;\lambda \right) \right]
\end{equation*}%
where

\begin{equation*}
\mathcal{P}_{3}\left( x-y;\lambda \right) =-\frac{1}{4\pi \lambda ^{4}}%
\underset{k=1}{\overset{3}{\sum }}\frac{\nu _{k}^{2}-A_{2}\nu _{k}}{\overset{%
3}{\underset{\underset{s\neq k}{s=1}}{\prod }}\left( \nu _{k}-\nu
_{s}\right) }\mathcal{K}_{0}\left( \frac{\lambda \left\vert x-y\right\vert }{%
\sqrt{-\nu _{k}}}\right) ,
\end{equation*}

\begin{equation*}
\mathcal{P}_{3}^{\ast }\left( x-y;\lambda \right) =-\frac{1}{4\pi \lambda
^{4}}\underset{k=1}{\overset{3}{\sum }}\frac{A_{1}\nu _{k}-A_{0}}{\overset{3}%
{\underset{\underset{s\neq k}{s=1}}{\prod }}\left( \nu _{k}-\nu
_{s}\right) }\mathcal{K}_{0}\left( \frac{\lambda \left\vert
x-y\right\vert }{\sqrt{-\nu _{k}}}\right) ,
\end{equation*}%
where $\mathcal{K}_{0}\left( z\right) $ is a second type Bessel
function of zero order.

By means of asymptotic and integral representations for the function $%
\mathcal{K}_{0}\left( z\right) $ and its derivatives (see [5]), we
first prove the necessary jump formulae for the potentials
$W_{k}\left( x,\lambda \right) \;\left( k=\overline{1,3}\right) $
and their derivatives for all  values of $\lambda \in R_{\delta
}$, where

\begin{equation*}
R_{\delta }=\left\{ \lambda :\left\vert \lambda \right\vert >R;-\frac{\pi }{4%
}+\delta \leq \arg \lambda <\frac{\pi }{4}\right\} .
\end{equation*}

By direct verification we prove that at all $\lambda \in R_{\delta }$ for
all potentials $w_{k}\left( x;\lambda \right) \;$

$\left( k=\overline{1,3}\right) $ it holds the estimation

\begin{equation*}
\left\vert \frac{\partial }{\partial x_{k}}\Delta ^{m}\mathcal{P}_{s}\left(
x-y;\lambda \right) \right\vert \leq \frac{Ce^{-\varepsilon \left\vert
\lambda \right\vert \left\vert x-y\right\vert }}{\left\vert \lambda
\right\vert ^{4-2m}\left\vert x-y\right\vert }\eqno(8)
\end{equation*}%
$\left( s=\overline{0,2},\;m=\overline{0,2},\;k=1,2\right) $ .

Putting (4) to the left hand sides of boundary conditions (2) and
taking into account the known jump formulae for the unknown
densities we obtain the following system of integral equations

\begin{equation*}
\mu \left( z,\lambda \right) =f\left( z,\lambda \right) +\underset{\tau }{%
\int }\mathcal{K}\left( z;y;\lambda \right) \mu \left( y;\lambda \right)
d\tau _{y},\ \eqno(9)
\end{equation*}%
where $\mu \left( z,\lambda \right) $ and $f\left( z;\lambda
\right) $are the columns of the functions composed of unknown
densities and boundary functions, and $\mathcal{K}\left(
x;y;\lambda \right) $ is a matrix of functions whose elements are
fundamental solution and particular solution and their
derivatives. By means of estimates (8) it is proved that for the
kernels $\mathcal{K}\left( x;y;\lambda \right) $ it holds an
estimation of type (8)  for all $\lambda \in R_{\delta }$.
Consequently, the system of integral equations is of Fredholm
property, so, we can solve it by the method of sequential
approximations and the solution is an analytic, bounded function
with respect to $\lambda $ in the domain $R_{\delta }$. So we
prove the

\textbf{Theorem:} \textit{Under conditions 1), 2), 3) problem (1), (2) has a
unique solution }$u\left( x,\lambda \right) $\textit{\ represented in the
form of the sum of potentials whose densities are the solutions of the
system of regular integral equations (9)}

\begin{center}
\textbf{References}
\end{center}

\bigskip

[1]. Rasulov M.L. \textit{Contour integral method.} "Nauka", M.:
1964

[2]. Rasulov M.L. \textit{Application of the contour integral
method.} "Nauka", M.: 1975.

[3]. Asadova O.H. \textit{On a plane problem for an elliptic
equation of the fourth order.} Dep. in AzNIINTI, 339, 1987, 12p.

[4]. Asadova O.H. \textit{On potentials of a sixth order
polyharmonic equation.} Dep. in VINITI, 1988, 24p.

[5]. Tikhonov A.N., Samarskii A.A. \textit{Equations of mathematical physics.%
} "Nauka", M.: 1966.

\end{document}